\documentclass[12pt]{article}
\title{Reverse mathematics of the finite downwards closed subsets of $\N^k$ ordered by inclusion and adjacent Ramsey for fixed dimension \footnote{This is a pre-peer-reviewed version of a paper which has been accepted for publication at Mathematical Logic Quarterly.} }
\author{Florian Pelupessy  \\
Mathematical Institute,
Tohoku University
}

\usepackage[english]{babel}
\usepackage{amsmath, amssymb}
\usepackage{url}
\usepackage{xcolor}
\usepackage{hyperref}
\usepackage{tikz}
\usepackage{multicol}
\usepackage{parskip}
\usepackage{makeidx}
  \makeindex
\usepackage{pdfcomment}
\usepackage{fullpage}

 \usepackage[initials]{amsrefs}

\BibSpec{article}{%
  +{}{\PrintAuthors}				{author} %
  +{,}{ \textit}					{title}
  +{,}{ }							{journal}
  +{}{ \textbf}					{volume}
  +{,}{ }							{note}
  +{}{ \parenthesize}				{date}
  +{:}{ }							{pages}
  +{.}{}							{transition}
  +{}{ }							{review}
  +{\\}{ \url }							{url}
}

\BibSpec{thesis}{%
  +{}{\PrintAuthors}				{author} %
  +{,}{ \textit}					{title}
  +{,}{ }							{journal}
  +{}{ \textbf}					{volume}
  +{,}{ }							{note}
  +{,}{ }							{date}
  +{:}{ }							{pages}
  +{.}{}							{transition}
  +{}{ }							{review}
  +{\\}{ \url }							{url}
}

\usepackage{fancyhdr}
\makeatletter
\usepackage[T1]{fontenc}
\usepackage{libertine}

\newtheorem{theorem}{Theorem}
\newtheorem{lemma}[theorem]{Lemma}
\newtheorem{definition}[theorem]{Definition}

\newtheorem{corollary}[theorem]{Corollary}

\newenvironment{remark}[1][Remark]{\begin{trivlist}
\item[\hskip \labelsep {\bfseries #1}]}{\end{trivlist}}

\newcommand{\qed}{
\begin{flushright}
$\Box$
\end{flushright}
}

\newcommand{\conc}{%
  \mathord{
    \mathchoice
    {\raisebox{1ex}{\scalebox{.7}{$\frown$}}}
    {\raisebox{1ex}{\scalebox{.7}{$\frown$}}}
    {\raisebox{.7ex}{\scalebox{.5}{$\frown$}}}
    {\raisebox{.7ex}{\scalebox{.5}{$\frown$}}}
  }
}

\newcommand{\RCA}{\mathrm{RCA}_0}
\newcommand{\RCAst}{\mathrm{RCA}_0^{\displaystyle{*}}}
\newcommand{\id}{\mathrm{id}}
\newcommand{\N}{\mathbb{N}}

\newcommand{\nin}{\not\in}

\begin{document}
\maketitle

\begin{abstract}
We show that the well partial orderedness of the finite downwards closed subsets of $\N^k$ ,ordered by inclusion, is equivalent to the well foundedness of the ordinal $\omega^{\omega^\omega}$. Since we use Friedman's adjacent Ramsey theorem for fixed dimensions in the upper bound, we also give a treatment of the reverse mathematical status of that theorem..
\end{abstract}


\section{Introduction}

In this note we prove the following theorem which was conjectured by Hatzikiriakou and Simpson in Remark 6.2 in \cite{HatzikiriakouSimpson}.


\begin{definition}
We order $k$-tuples coordinatewise.
\end{definition}

\begin{theorem}\label{thm:main}
$\RCA$ proves that the following are equivalent:
  \begin{enumerate}
  \item $\omega^{\omega^\omega}$ is well founded,
  \item For every $k$: the finite downwards closed subsets of $\N^k$, ordered by inclusion, are a well partial order.
  \end{enumerate}
\end{theorem}

The case of $k=2$  was shown by Hatzikiriakou and Simpson to be equivalent to the well-foundedness of $\omega^\omega$ \cite{HatzikiriakouSimpson}. As remarked in that paper, there is an order-preserving one-to-one correspondence between $k$-dimensional partitions from Chapter~11 from \cite{andrews} and the  downwards closed finite sets in $\N^{k+1}$. We will examine this in more detail at the end of section \ref{section:downwardsclosed}. This note is also part of the attention, in reverse mathematics, for the strength of the well foundedness of the ordinals $\omega^\omega$ and $\omega^{\omega^\omega}$ (See, e.g.: \cite{simpson2015, kreuzeryokoyama, HatzikiriakouSimpson}). Research on these levels goes back to Simpson's work on the Hilbert and Robson basis theorems, in \cite{simpson1988},  or even further to Goodstein's work on his sequences, in \cite{goodstein}. 

We assume basic familiarity with reverse mathematics in $\RCA$ (II.1-II.3 in \cite{sosoa}) and treatment of ordinals less or equal to $\omega^{\omega^\omega}$ and their Cantor Normal Forms (See, e.g. Definition~2.3 in \cite{simpson1988} or Section~II.3(a) in \cite{hajekpudlak}). The remainder of this note is divided in two sections: one on finite downwards closed subsets and the other on adjacent Ramsey. In the latter we  treat the supporting  upper bound used in the other, which may also be of interest in its own right. 

\section{Finite downwards closed subsets}\label{section:downwardsclosed}
\subsection{Equivalence}
\begin{definition}
Given partial order $(X, \leq)$, we call a sequence $x_0, x_1, \dots$ of elements from $X$ bad if for all $i<j$ we have $x_i \not\leq x_j$. 
\end{definition}

\begin{definition}
A partial order is a well partial order (w.p.o.) if every bad sequence in the order is finite. 
\end{definition}

\begin{definition}
A partial order is well founded if every strictly descending sequence in the order is finite.
\end{definition}

We will use the following principle from Friedman \cite{friedman2010} for the upper bound:
\begin{definition}[Adjacent Ramsey for pairs]
For every function $C\colon \N^2 \rightarrow N^r$ there exist $a<b<c$ with $C(a,b) \leq C(b,c)$.
\end{definition}

\begin{theorem}\label{thm:aramseyequiv}
$\RCA$ proves that the following are equivalent:
  \begin{enumerate}
  \item $\omega^{\omega^\omega}$ is well-founded,
  \item the adjacent Ramsey theorem for pairs.
  \end{enumerate}
\end{theorem}
\emph{Proof:} See next Section~\ref{section:adjacent}.
\qed

\emph{Proof of Theorem~\ref{thm:main} (1) $\rightarrow$ (2):} Take, for a contradiction, an infinite bad sequence $G_0 , G_1, G_3, \dots$ with:
\[
G_i=\{ m_{i, 0} \dots , m_{i, n_i}\}
\]
Define:
\[
C(i, j) = m_{i,l}
\]
where $l\leq n_i$ is the smallest such that $\forall p\leq n_j .m_{i,l} \not\leq m_{j,p}$. By adjacent Ramsey there exist $a<b<c$ such that $C(a,b) \leq C(b,c)$, contradiction.
\qed

\begin{definition}
A downwards closed subset $X$ is generated by $G$  if 
\[
X=\{ m \in \N^k : \exists m'\in G.m \leq m' \}.
\]
\end{definition}

Every finitely generated set is also finite (upper bound given by the generators).

For finite sets we say $G\leq H$ if $X\subseteq Y$, where $X$ and $Y$ are the respective generated sets.

Notice that $G\not\leq H$ if and only if there exists $m \in G$ with $\forall m' \in H . m \not\leq m'$. 

\emph{Proof of Theorem~\ref{thm:main} (2) $\rightarrow$ (1):} For $\beta = \omega^k \cdot b_0 + \dots + \omega^0 \cdot b_k < \omega^{k+1}$, take $h(\beta)= (b_0, \dots , b_k) \in \N^{k+1}$. We have the following property: $h(\beta) \leq h(\beta') \rightarrow \beta \leq \beta'$.

For $\alpha=_{\mathrm{CNF}}\omega^{\beta_0} \cdot a_0 + \dots + \omega^{\beta_n}\cdot a_n< \omega^{\omega^{k+1}}$, define:
\[
f(\alpha)=\{ (i, a_i)\conc h(\beta_i)  : i \leq n \}.
\]
Notice that $f(\alpha)$ is an antichain in $\N^{k+3}$.

Assume, for a contradiction, that $\omega^{\omega^{k+1}} > \alpha_0 > \alpha_1 >\dots$ is an infinite sequence and let $i<j$ be such that $f(\alpha_i) \leq f(\alpha_j)$ by well-partial-orderedness. Denote:
\[
\alpha_i=_{\mathrm{CNF}}\omega^{\beta_{i,0}} \cdot a_{i, 0} + \dots + \omega^{\beta_{i, n_i}}\cdot a_{i,n_i},
\]
\[
\alpha_j=_{\mathrm{CNF}}\omega^{\beta_{j,0}} \cdot a_{j, 0} + \dots + \omega^{\beta_{j, n_j}}\cdot a_{j,n_j}.
\]
Let $l$ be the smallest such that $(l, a_{i,l})\conc h(\beta_{i,l}) \not\leq (l, a_{j,l})\conc h(\beta_{j,l})$, such $l$ exists because otherwise $\alpha_i \leq \alpha_j$.

Let $q>l$ be the smallest such that $(l, a_{i,l})\conc h(\beta_{i,l}) \leq (q, a_{j,q})\conc h(\beta_{j,q})$, such $q$ exists because of $f(\alpha_i) \leq f(\alpha_j)$.

By the properties of the Cantor Normal Forms, we have the following for all $p \geq l$:
\[
\omega^{\beta_{j,l}} >  \omega^{\beta_{j,q}} \geq \omega^{\beta_{i,l}} \geq \omega^{\beta_{i,p}}.
\]
Hence, by $\omega^{\beta_{j,l}}$ being closed under ordinal addition, $\alpha_i \leq \alpha_j$, contradiction.
\qed

\subsection{Higher dimensional partitions}
We turn our attention to the consequences of the previous section for the $k$-dimensional partitions from Chapter~11 of \cite{andrews}. 

\begin{definition}
A $k$-dimensional partition $N$ of $n$ is a term
\[
\sum_{(i_1, \dots , i_k)\in A} n_{i_1, \dots , i_k},
\]
with the following properties:
\begin{enumerate}  
   \item $A$ is downwards closed,
  \item the $n$'s are strictly positive integers, occuring in the expression in lexicographic order of the subscipts,
  \item $n$ is the value of the term, using the canonical interpretation of sums,
  \item if $i_1 \leq j_1$, $\dots$, $i_k \leq j_k$ then $n_{i_1, \dots , i_k}\geq n_{j_1, \dots , j_k}$.
\end{enumerate}
We denote the value of $N$ with $v(N)$.
\end{definition}

We generalise the ordering as given for the one dimensional case in \cite{HatzikiriakouSimpson}.
\begin{definition}
Given 
\[
N=\sum_{(i_1, \dots , i_k)\in A} n_{i_1, \dots , i_k}, M=\sum_{(i_1, \dots , i_k)\in B} m_{i_1, \dots , i_k},
\]
we write $N\leq_t M$ if $n_{i_1, \dots , i_k} \leq m_{i_1, \dots , i_k}$ for all $(i_1, \dots , i_k)\in A$, where $m_{i_1, \dots , i_k}$ is read as $0$ whenever $(i_1, \dots , i_k)\nin B$.
\end{definition}

Notice that if $N\leq_t M$, then 
\[
\{(n_{i_1, \dots , i_k}-1, i_1, \dots , i_k): (i_1, \dots , i_k) \in A \}\leq\{(m_{i_1, \dots , i_k}-1, i_1, \dots , i_k): (i_1, \dots , i_k) \in B\}.
\]
Inversely, if $X,Y \in D_{k+1}$ and $X\subseteq Y$, then 
\[
\sum_{(i_1, \dots , i_k) \in A} n_{i_1, \dots , i_k} \leq_t \sum_{(i_1, \dots , i_k)\in B} m_{i_1, \dots , i_k},
\]
where $n_{i_1, \dots , i_k}=\max \{n+1: (n, i_1, \dots , i_k) \in X\}$, $A=\{(i_1, \dots , i_k): \exists i (i, i_1, \dots , i_k)\in X\}$  and $m$'s and $B$ taken similarly from $Y$. Hence, we can generalise the one dimensional partitions:

\begin{corollary} \label{corollary:partitionsequiv}
$\RCA$ proves that the following are equivalent:
  \begin{enumerate}
  \item $\omega^{\omega^\omega}$ is well founded,
  \item For every $k \in \N$, the $k$-dimensional partitions, ordered by $\leq_t$, are a well partial order.
  \end{enumerate}
\end{corollary}

\begin{remark}
Given that our upper bound proof for adjacent Ramsey is based on the one for the first order variant, we can simply observe the upper bound for the Friedman-style miniaturisation of the well orderedness of partitions in the following manner:

\begin{corollary} \label{corollary:miniprov}
Given $k$, the following is provable in $\mathrm{I}\Sigma_2$: For every $l \in \N$ there exists $R$ such that for every sequence $N_0, \dots , N_R$ of $k$-dimensional partitions, with $v(N_i) \leq l+i$, there are $i<j\leq R$ with $N_i \leq_t N_j$.
\end{corollary}

Furthermore, we have the following:

\begin{corollary} 
Given $k>0$ standard and $f\colon \N\rightarrow\N$, there exists $g\colon \N\rightarrow\N$, multiply recursive in $f$ and $l$, such that every bad sequence of $k$-dimensional partitions $N_0, N_1, \dots$ with $v(n_i) \leq f(l+i)$ has maximum length $g(l)$.
\end{corollary}

We expect that for unrestricted $k>0$, already for $f=\id$, there is no such function which is multiply recursive in $k$ and $l$. Furthermore we expect Corollary~\ref{corollary:miniprov} not to hold for the statement with unrestricted dimension.

\end{remark}

\section{Adjacent Ramsey}\label{section:adjacent}
In this section we prove the upper bound for Theorem~\ref{thm:aramseyequiv}. On a side note, if one does not need the tight upper bound, it is possible to easily prove adjacent Ramsey directly from Ramsey's Theorem. Friedman used this fact in his proof of the upper bound for adjacent Ramsey with unrestricted dimensions \cite{friedman2010}, which makes his proof not suitable for use in the case of a fixed dimension.

Since this requires no extra effort, we will be treating the general case for arbitrary dimension $d+1$. The proof for the upper bound is a simple assembly of adaptations of the proofs for existing first order, finitary results. We start with a few definitions:

\begin{definition} ~{}
\begin{enumerate} 
  \item $\omega_1=\omega$, $\omega_{n+1}=\omega^{\omega_n}$, $\omega_1 (\alpha)=\omega^\alpha$, $\omega_{n+1}(\alpha)=\omega^{\omega_n(\alpha)}$,
  \item We use terminology from Ramsey theory: $[X]^d$ is the set of $d$-element subsets of $X$, $[a,R]^d=[\{a, \dots , R \}]^d$, and we identify any $c \in \N$ with $\{0, \dots , c-1\}$.
  \item Given a colouring $C\colon [X]^d \rightarrow c$, we call $H\subseteq X$ homogeneous for $C$, or $C$-homogeneous, if $C$ is constant on $[H]^d$.
  \item Given a colouring $C\colon [X]^{d+1} \rightarrow c$, we call $H=\{h_0 < \dots < h_n \}\subseteq X$ adjacent homogeneous for $C$, or $C$-adjacent-homogeneous, if $C(h_i , \dots , h_{i+d}) = C(h_{i+1} , \dots , h_{i+d+1})$ for all $i < n-d$. 
  \item A colouring $C\colon \{0, \dots R\}^d \rightarrow \N^r$ is $f$-limited if 
  \[
  \max C(x_1, \dots x_d) \leq f(\max \{x_1, \dots , x_d)\}).
  \]
\end{enumerate}
\end{definition}

\begin{theorem} \label{thm:equivalences}
The following is provable in $\RCA$:  for every $d$, the following are equivalent:
\begin{enumerate}
  \item $\omega_{d+2}$ is well founded, 
  \item  the parametrised Paris--Harrington principle in dimension $d+2$: for any $f\colon \N \rightarrow \N$, $a,c \in N$ there exists $R$ such that for every $C \colon [a,R]^{d+2} \rightarrow c$ there is a $C$-homogenous $H \subseteq [a,R]$ of size $> f(\min H)$,
  \item the parametrised adjacent Paris--Harrington principle in dimension $d+2$: for any $f\colon \N \rightarrow \N$, $a,c \in N$ there exists $R$ such that for every $C \colon [a,R]^{d+2} \rightarrow c$ there is a $C$-adjacent-homogeneous $H \subseteq [a,R]$ of size $> f(\min H)$,
  \item the parametrised strong adjacent Paris--Harrington principle in dimension $d+2$: for any $f\colon \N \rightarrow \N$, $a,c,k \in N$ there exists $R$ such that for every $C \colon [a,R]^{d+2} \rightarrow c$ there is a $C$-adjacent-homogeneous $H =\{ h_0< \dots < h_{|H|-1}\} \subseteq [a,R]$ of size $> f(\min h_k)$,
  \item the parametrised finite adjacent Ramsey theorem: for any $f\colon \N \rightarrow \N$, $r \in \N$ there exists $R$ such that for every $f$-limited $C\colon \{0,\dots, R\}^{d+1} \rightarrow \N^r$ there are $x_1 < \dots <x_{d+2}$ with $C(x_1, \dots , x_{d+1}) \leq C(x_2, \dots , x_{d+2})$,
  \item adjacent Ramsey in dimension $d+1$: for every $C\colon \N^{d+1} \rightarrow \N^r$ there are $x_1 < \dots < x_{d+2}$ with $C(x_1, \dots , x_{d+1}) \leq C(x_2, \dots , x_{d+2})$.
\end{enumerate}
\end{theorem}

The first order variant of (1)$\rightarrow$(2) is due to Ketonen and Solovay \cite{ketonensolovay}. The first order variant of (2)$\rightarrow$(3)$\rightarrow$(4)$\rightarrow$(5) is due to Friedman \cite{friedmanpelupessy}. (5)$\rightarrow$(6) is self evident. (6)$\rightarrow(1)$ is a modification of Friedman's treatment in \cite{friedman2010} for the unrestricted dimensions and $\varepsilon_0$. The first order variant of (1) is the totality of the function $H_{\omega_{d+2}}$ from the Hardy hierarchy, for the other items the first order variant is obtained by restricting $f$ to just the identity function. 

We first use the concept of $\alpha$-largeness from \cite{ketonensolovay}.

\begin{definition}
A finite set $A=\{a_0 < \dots < a_b \}$ is $\alpha$-large if $\alpha[a_0] \dots [a_b]=0$, where $\alpha [.]$ denotes the canonical fundamental sequences for ordinals below $\varepsilon_0$.
\end{definition}

The key ingredient for (1)$\rightarrow$(2) is Theorem~6.7 from \cite{ketonensolovay}. By a straightforward verification, the proof of this theorem in \cite{ketonensolovay} is within $\RCA + $ ``$\omega_{d+2}$ is well founded'':

\begin{theorem}[Ketonen--Solovay]\label{thm:ketonensolovay}
If $A>3$ is $\omega_{d+1}(c+5)$-large, then for any $D\colon [A]^{d+2} \rightarrow c$ there exists $H \subseteq A$ of size $> \min H$ such that $D$ is constant on $[H]^{d+2}$.
\end{theorem}

\textbf{Note:} that the Ketonen--Solovay proof has many applications of transfinite induction which are all consequences of the well-foundedness of the appropriate ordinal. This is sufficient for our purpose. The interested reader can find in \cite{pelupessyketonensolovay} a description of how to remove all instances of transfinite induction.

\begin{lemma}\label{lemma:simple}
$\RCA$ proves the following:  if $\omega_{d+2}$ is well-founded, then  for every strictly increasing  $f\colon \N \rightarrow \N$, $a \in \N$, $\alpha < \omega_{d+2}$ there exists $\alpha$-large set $\{ f(a) , f(a+1) ,\dots , f(b) \}$.
\end{lemma}
\emph{Proof of the lemma:} Take the following descending sequence of ordinals: $\alpha_0=\alpha$ and:
\[
\alpha_{i+1}=\alpha_i [f(i)].
\]
By well-foundedness of $\omega_{d+2}$ this sequence reaches zero, delivering the desired $\alpha$-large set.
\qed

\begin{lemma}
(1)$\rightarrow$(2)
\end{lemma}

\emph{Proof:} Assume, without loss of generality, that $f$ is strictly increasing. Take $\omega_{d+1}(c+5)$-large $A=\{f(a) , \dots , f(b) \}$ from Lemma~\ref{lemma:simple}, $R=b$, and arbitrary $C\colon [a,R]^{d+2} \rightarrow c$. Define $D(x_1,\dots, x_{d+2})= D(f^{-1} (x_1), \dots , f^{-1} (x_{d+2}))$ on $A$.

By Theorem~\ref{thm:ketonensolovay}, there exists $\bar{H} \subseteq A$ of size $> \min \bar{H}$ such that $D$ is constant on $[\bar{H}]^{d+2}$. Then $H=\{ f^{-1} (h) : h \in \bar{H}\}$ is the desired subset of $[a,R]$.
\qed

\begin{lemma}
(2)$\rightarrow$(3)$\rightarrow$(4)
\end{lemma}
\emph{Proof:} (2)$\rightarrow$(3) is trivial. Assume, without loss of generality, that $f$ is strictly increasing and that $k>0$, $a> d+k+2$. Take $R$ from the adjacent Paris--Harrington principle with $f$ , but with codomain $2c$.  Given $C\colon [a, R]^{d+2} \rightarrow c$, define the following colouring:

$D_1(x_1, \dots x_d)=1$ if there exist $z_0 < \dots < z_{k-1} < x_1$ such that $\{z_0 , \dots , z_k , x_1, \dots , x_d\}$ is $C$-adjacent-homogeneous, $0$ otherwise. 

Obtain $D$ by combining $D_1$ and $C$ into a single function with codomain $2c$. Observe that for any $D$-adjacent-homogeneous $H$ of size $>f(\min H)$, by definition of $D_1$, there exist $z_0 < \dots < z_{k-1}$ such that $\{z_0, \dots , z_{k-1}\} \cup H$ is the desired $C$-adjacent-homogenous set. 
\qed

\begin{lemma}
(4)$\rightarrow$(5)
\end{lemma}
\emph{Proof:} Given $r$, assume without loss of generality, that $f$ is strictly increasing, take $a=d+4$, $k=d$ and $R$ from the strong adjacent Ramsey principle with codomain $r+1$. Let $C\colon \{0, \dots , R\}^{d+1} \rightarrow \N^r$ be $f$-limited. Take:  
 \[
D(x_1, \dots , x_{d+2})=\left\{ 
\begin{array}{ll}
0 & \textrm{if $C(x_1-a, \dots , x_{d+1}-a) \leq C(x_2-a, \dots , x_{d+2}-a)$}, \\
i  & \textrm{otherwise},
\end{array}
\right.
\]
where $i$ is the least such that 
\[
(C(x_1-a, \dots , x_{d+1}-a))_i > (C(x_2-a, \dots , x_{d+2}-a))_i.
\]
By the choice of $R$, there is $D$-adjacent-homogenous $H=\{ h_0 < \dots < h_{f(h_d)} \}$. If $D(h_0, \dots , h_{d})\neq 0$ we obtain a strictly descending sequence starting with $m\leq f(h_{d}-a)\leq f(h_{d})-a$ of length $f(h_d)-d$, which  is impossible. Hence:
\[
C(h_0-a, \dots ,h_{d}-a) \leq C(h_2-a, \dots, h_{d+1}-a). 
\]
\qed
\begin{lemma}
(5)$\rightarrow$(6)
\end{lemma}
Given $C\colon \N^{d+1} \rightarrow \N^r$, take $f(x)= \max_{\bar{y} \in \{ 0, \dots , x\}^{d+1}} C(\bar{y})$ to obtain the desired $x_1 < \dots < x_{d+2}\leq R$ from (5). 
\qed
\begin{lemma}
(6)$\rightarrow$(1)
\end{lemma}
See Definitions and Lemmas 1.8-1.11 and the first three lines of the proof of Theorem~2.1 from \cite{friedmanpelupessy}, but with an arbitrary sequence of ordinals below $\omega_{d+1}(l)$. 
\qed

\end{document}